# Finding Your Way: Shortest Paths on Networks

By Teresa Rexin and Mason A. Porter

## Abstract

Traveling to different destinations is a big part of our lives. We visit a variety of locations both during our daily lives and when we're on vacation. How can we find the best way to navigate from one place to another? Perhaps we can test all of the different ways of traveling between two places, but another method is to use mathematics and computation to find a shortest path. We discuss how to construct a shortest path and introduce Dijkstra's algorithm to minimize the total cost of a path, where the cost may be the travel distance, travel time, or some other measurement. We also discuss how to use shortest paths in the real world to save time and increase traveling efficiency.

## What is a path?

Every day, we make decisions about which routes we use to travel between different places. In your house, you may travel from your bedroom to your kitchen. Outside your house, you may travel from your home to school. Suppose that we have a **network** of places that are connected to each other by streets, walkways, and other things that we use to travel between locations. Each of these locations is called a **node**, and the streets and walkways are called **edges**. The **neighbors** of a node are the nodes to which it is connected by an edge. A **path** is a sequence of edges between an origin node and a destination node [1, 2].

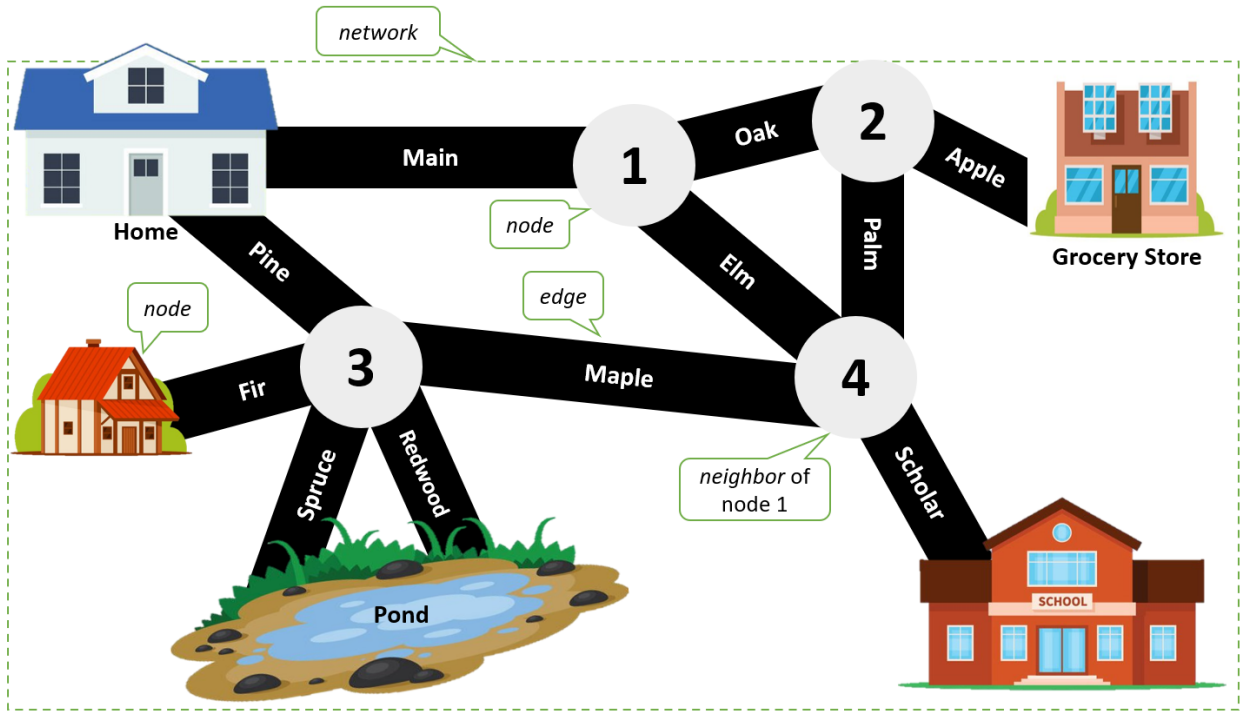

**Figure 1:** A city map inspired by Clipart.Email Map Street Clipart (from https://www.clipart.email/clipart/map-street-clipart-390373.html). Locations (in other words, nodes) 1–4 occur at the intersections between each street (that is, edge) between them. Each of the locations (the blue and red houses, the pond, the school, and the grocery store) is also a node.

*Activity 1:* In Figure 1, trace a path from the blue house to the school on the map. Which streets (that is, edges) do you take in the picture? Which edges do you take on your path from your house to school in real life?

# Shortest paths

In mathematics, one often studies the lengths of paths and tries to construct short paths. A **shortest path** is a path between two nodes that has the fewest edges if the **cost** of traveling along each edge is the same (for example, if each edge is a street of the same length). More generally, a shortest path from an origin node to a destination node is a path that has the smallest sum of edge costs of all of its edges, among all paths with the same origin and same destination [1]. A cost may measure distance, time, or something else. For example, in the city map in Figure 1, a shortest path from home to school may be one that takes the least amount of time among the possible paths. There can be more than one shortest path between two nodes in a network, as multiple paths may have the same minimum cost. That is why we refer to "a" shortest path between two nodes (even though it sounds weird) rather than "the" shortest path between them.

You probably already think about shortest paths in your daily life when you're going to different places. In our bedroom-to-kitchen example, it wouldn't make much sense to walk from your bedroom, then to the laundry room, then outside to your backyard, and finally to your kitchen if you only want to travel from your bedroom to your kitchen. It would be much faster to walk directly from your bedroom to your kitchen without stopping in the laundry room and your backyard first (unless perhaps you also have chores to do in those places).

In a trip with locations that are nearby, there are few enough street intersections (in other words, nodes) and you may be able to try out a large number of the different paths to find a shortest path. But if the locations are farther apart — say, your home, your school, and a toy store in a different city — then finding a shortest path is very difficult to do. How do navigation tools like Google Maps determine the best way to reach a destination? One way is to solve the **shortest-path problem**, which is the problem of finding a path between two nodes in a way that minimizes the sum of the costs of the edges in the path [4].

In mathematics, we often label nodes by using numbers (like the intersections in Figure 1) or by using letters (like in Figure 2). For simplicity, we also suppose that everything is two-dimensional (like a drawing on a piece of paper), so we'll measure distance the way we would between two spots on the floor in your house, rather than worrying about things like height or the curvature of the Earth. In the network in Figure 2, if we want to find a shortest path from node A to node F, we should choose the edges with the lowest costs. For example, instead of choosing the edge with cost 4 from node A to node B, we choose the edge with cost 2 from node A to node C. Choosing the edges with the lowest costs to find a shortest path is one of the key ideas in *Dijkstra's algorithm*[1] [5].

---

[1] When pronouncing the name Dijkstra, note that the 'j' is silent.

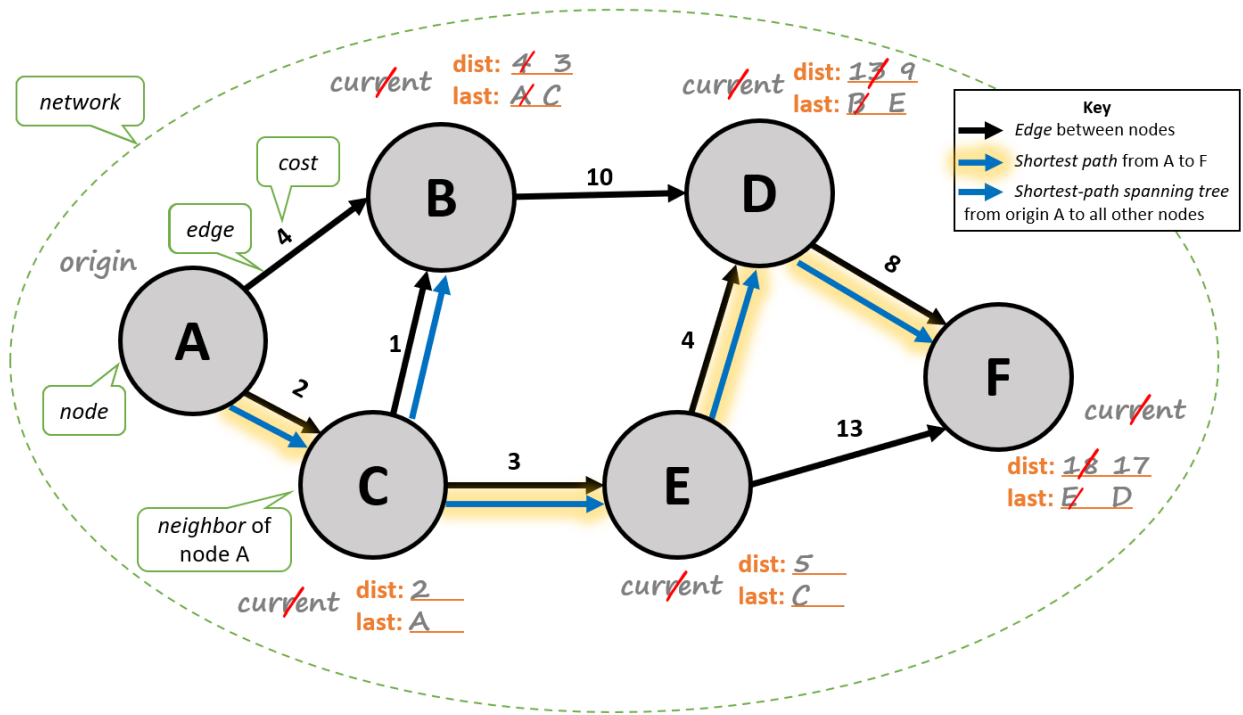

**Figure 2:** In this network, following the highlighted blue arrows shows us the shortest path from node A to node F. We use numbers to display the costs of the edges. (The edge lengths are not drawn to scale.) We use blue arrows to show the shortest-path spanning tree with A as the origin node. Notice that the shortest path from A to F is also part of the shortest-path spanning tree. We work through this example in the section called "Dijkstra's algorithm". Inspired by [4] and [5], the shorthand "dist" indicates the total distance from the origin node "origin" to a particular node and "last" indicates the last node that one passes through to reach a particular destination node from "origin".

# Dijkstra's algorithm

An **algorithm** is a precise set of steps to follow to solve a problem, such as the shortest-path problem [1]. **Dijkstra's algorithm** [5] is a famous shortest-path algorithm; it is named after its inventor Edsger Dijkstra [6], who was a Dutch computer scientist. One adaptation of Dijkstra's algorithm is to methodically create a **shortest-path spanning tree** to find shortest paths from an origin node to each other node in a network by calculating the distances one node at a time. In this algorithm, whenever we find a shorter path to a node through a neighboring node, we update the distance. We use distance for concreteness, but we can use Dijkstra's algorithm for any type of cost.

We now present the algorithm to create a shortest-path spanning tree for a connected network. (We adapted it from the description in [5] and are using Figure 2 as an example. See the video explanation at https://drive.google.com/file/d/1sQZHh0hQE6WBeXVuCv-i3rBVJ4ERyFXK/view to follow along.) We proceed as follows:

1. We shade in the origin node (labeled "origin"). For each of its neighbors, we set the initial value of "dist" to be the distance from origin to each neighbor of origin and the initial value of "last" to be the origin node. Node A is the origin node in Figure 2.
2. We identify the unshaded node with the lowest "dist" value (excluding blanks) and label this as our "current" node. *(Example: If we begin with node A as the origin, then the current node is node C. The reason is that nodes B and C are the only neighbors of A and the dist of 2 from A to C is less than the dist of 4 from A to B in Figure 2.)* If there is a tie, we choose any of the nodes with the smallest "dist" value.
3. We do the following steps for each unshaded neighbor of current:
   a. We add current's dist to the cost of the edge from current to the neighbor.
   b. If dist from Step 3a is smaller than the neighbor's dist (or if the neighbor's dist isn't written down yet), we update the neighbor's dist to the dist that we calculated in Step 3a and set the neighbor's "last" to be the letter of the current node.
4. After we complete Step 3 for all unshaded neighbors of current, we shade in current and cross out the label "current".
5. If all nodes are shaded, we go to Step 6. Otherwise, we go back to Step 2.
6. We highlight the edge between each node and its "last" node to reveal a shortest-path spanning tree from the origin.

*Activity 2*: Now it's your turn! Use Dijkstra's algorithm to find a shortest-path spanning tree from the origin to each other node in the network in Figure 3. (We have completed Step 1 for you.)

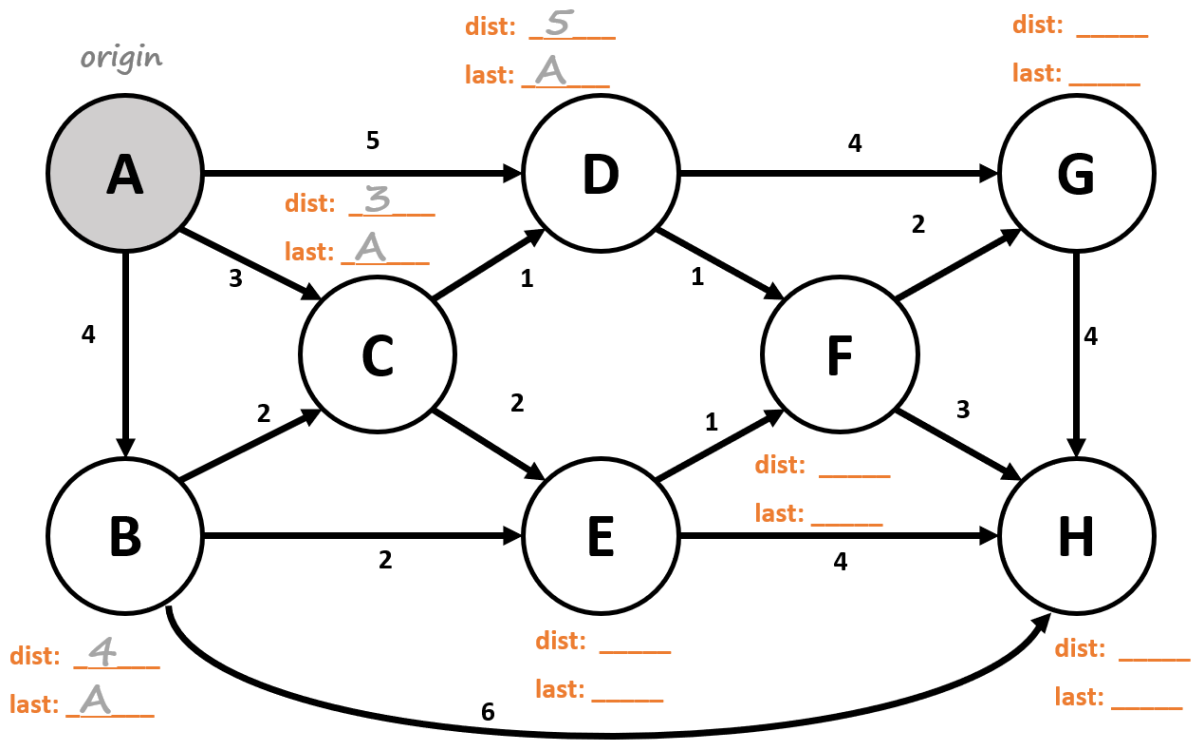

**Figure 3:** A network (inspired by [4]) for you to practice finding a shortest-path spanning tree. You can download a printable version of Figure 3 from https://drive.google.com/file/d/1rNONK-cmy4gq_aCJRnAYpe9Y2HSeq2A1/view.

# Applications

Using Dijkstra's algorithm, we can find a shortest path from an origin node to any other node in a network. If you think of your house as the origin node and your destination as a different node in a network, you can determine the fastest route from your home to any place that you want to go.

Suppose that you want to visit several places before returning home. How do you find the best way to visit all of the destinations while minimizing expenses, such as gas, hotels, and time? More abstractly, how can we find a shortest path that passes through all of the nodes in a network and returns to the starting node? This problem, which is an extension of the shortest-path problem, is known as the "**Traveling Salesperson Problem**"[2].

Finding shortest paths is also important for solving problems in many different types of networks. Shortest paths can improve the efficiency of city planning. For example, civil engineers can represent a city as a network and determine the best places to build different structures, such as roads to reduce traffic congestion and irrigation pipes to distribute the water supply to the population on demand [2]. Finding shortest paths also enables the transfer of data from one computer to another at high speeds, allowing massive amounts of information to travel in seconds [1, 2].

There are also many examples of shortest paths in communication and social networks. We briefly discuss two examples.

In our first example, suppose that each person is a node and each edge represents a friendship. You can find how to connect to a person outside of your friendship groups through other people's connections. This was explored in experiments by Stanley Milgram and his collaborators. In their experiments, they illustrated that the paths of connections (such as friendships) between two random people in the United States are shorter than one might think and also — astoundingly! — that people tend to be very good at navigating these paths with very little information [7]. The shortness of these paths is known as the "small-world phenomenon", and the small path lengths (with fewer than six steps, on average, between an origin person and a destination person in a path) also inspired the term "six degrees of separation" [1].

---

[2] See Cook, W. 2018. "Information, computation, optimization: Connecting the dots in the traveling salesman problem", which is at https://www.youtube.com/watch?v=q8nQTNvCrjE&t=35s.

Our second example relates to current events. During the COVID-19 pandemic, finding shortest paths has been valuable for limiting exposure to others when performing essential tasks. When moving in supermarkets, for instance, it is helpful to find a shortest path to pick up your groceries while avoiding contact with others through physical distancing [8, 9, 10].

# Conclusions

Shortest paths are an interesting and important idea when thinking about traveling in your neighborhood and in the world. They have numerous applications in networks of all sorts and can help solve a variety of real-world problems. From planning a family vacation to exploring how our world is connected, the study of shortest paths on networks is incredibly important and forms the basis for more complex investigations.

# Glossary

**Algorithm.** A set of precise steps to follow to solve a problem. An example of an algorithm is Dijkstra's algorithm in Activity 2.

**Cost.** A measure of how much effort it takes to travel along an edge in a network. In real life, a cost may measure distance, time, or something else.

**Edge.** An object that connects two nodes to each other. For example, when going from your house to school in Figure 1, each street is an edge.

**Neighbor.** Given a node, the other nodes to which it is connected by an edge are the neighbors of that node.

**Network.** A collection of objects (in other words, nodes) and the connections (in other words, edges) between those nodes.

**Node.** The objects in a network that are connected to other objects. For example, when going from your house to school in Figure 1, each location and street intersection is a node.

**Path.** A sequence of edges from an origin node to a destination node.

**Shortest path.** A path from an origin node to a destination node that has the lowest total cost among all paths from the origin to the destination.

**Shortest-path problem.** The problem of constructing a path between two nodes in a way that minimizes the sum of the costs of the edges in the path.

**Shortest-path spanning tree.** A subset of a connected network that indicates the shortest path from a specified origin node to any other node in the network.

**Traveling salesperson problem.** The problem of constructing a shortest path that goes through all nodes of a network and then returns to the origin node.

# Answer Key

*Activity 1:* One possible path is (Main, Elm, Scholar). Another possible path is (Main, Oak, Palm, Scholar). A third possible path is (Pine, Maple, Scholar).

*Activity 2:* Here is one example of a completed shortest-path spanning tree on the network in Figure 3. Follow along with this video explanation: https://drive.google.com/file/d/1b5rdzMDLBfvsnyTaCu-uEJ_4za5yuw4n/view.

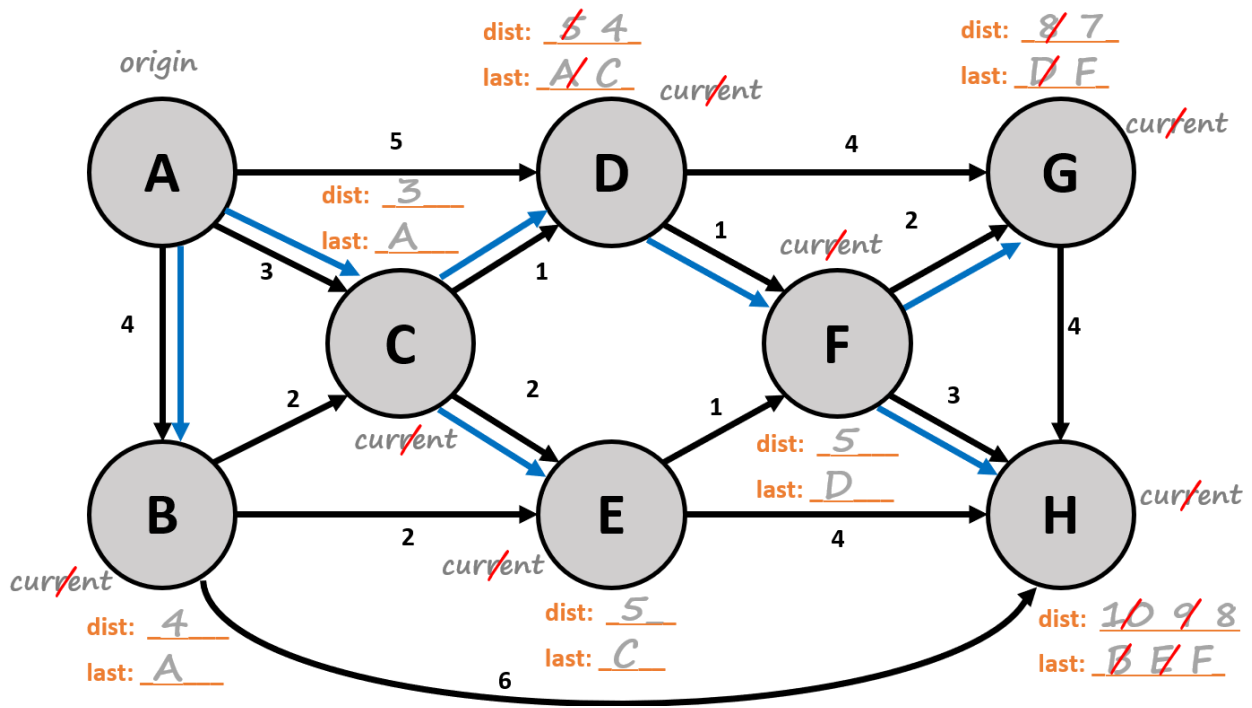

# References

[1] Newman, M. E. J. 2018. *Networks*, 2nd Edn. Oxford: Oxford University Press.
[2] NetSciEd. (Eds). 2015. Network literacy: Essential concepts and core ideas. Available online at http://tinyurl.com/networkliteracy
[4] Wikipedia. 2020. Shortest path problem. Available online at: https://en.wikipedia.org/wiki/Shortest_path_problem (Accessed 20 Aug, 2020).
[5] Code.org. 2020. U2L07 Activity guide — Dijkstra's shortest path algorithm. Available online at

https://docs.google.com/document/d/15N7aHAoWG1_9VIcDHNZRygzFK0hle-EHlmHu0PZI8D4/view (Accessed 20 Aug, 2020).
[6] Wikipedia. 2020. Edsger W. Dijkstra. Available online at https://en.wikipedia.org/wiki/Edsger_W._Dijkstra (Accessed 20 Aug, 2020).
[7] Milgram, S. 1967. The small-world problem. *Psychology Today*, **1**(1), 61–67.
[8] Ying, F., Wallis, A. O. G., Beguerisse-Díaz, M., Porter, M. A., and Howison, S. D. 2019. Customer mobility and congestion in supermarkets. *Phys. Rev. E* **100**, 062304.
[9] Brooks, H. Z., Kanjanasaratool, U., Kureh, Y. H., and Porter, M. A. 2020. Disease detectives: Using mathematics to forecast the spread of infectious diseases. *Frontiers for Young Minds* **8**, 577741. https://kids.frontiersin.org/article/10.3389/frym.2020.577741
[10] Ying, F. and O'Clery, N. 2020. Modelling COVID-19 transmission in supermarkets using an agent-based model, https://arxiv.org/abs/2010.07868.

# Conflict of interest statement

The authors declare that the research was conducted in the absence of any commercial or financial relationships that could be construed as a potential conflict of interest.

# Acknowledgements

We are grateful to our young readers — Nia Chiou, Taryn Chiou, Zoë Chiou, Tycho Elling, and Sage Hansen — for their many helpful comments. We also thank their family members — Lyndie Chiou, Christina Chow, Tim Elling, and Sterling Hansen — for putting us in touch with them and soliciting their feedback. We also thank Lyndie Chiou, Michelle Lee, Thomas Rexin, and Akrati Saxena for helpful comments. We also thank our young reviewers and their mentors for their many excellent suggestions. MAP acknowledges support from the National Science Foundation (grant number 1922952) through the Algorithms for Threat Detection (ATD) program.

# Author Biographies

**Mason A. Porter** is a professor in the Department of Mathematics at UCLA. He was born in Los Angeles, California, and he is excited to be a professor in his hometown. In addition to studying networks and other topics in mathematics and its applications, Mason is a big fan of games of all kinds, fantasy, baseball (Go Dodgers!), the 1980s, and other delightful things. Mason used to be a professor at University of Oxford, where he did actually wear robes on occasion (like in the *Harry Potter* series). Mason's most common shortest paths occur between his apartment and places to get good coffee.

**Teresa Rexin** was raised in Sacramento County, California. Teresa just finished her undergraduate degree at UCLA in applied mathematics and statistics and will be starting her

Master's in Statistics at UCSD this fall. She is interested in applying mathematics to solve challenging everyday problems and help improve people's quality of life. In her free time, she enjoys volunteering in her community, working out at the gym, and spending time with friends and family. Although it was not a shortest path back to her dorm from the UCLA campus, Teresa misses walking down UCLA's Bruin Walk and watching the sunset after a long week of classes.